\documentclass[11pt]{article}

\usepackage{amsfonts}
\usepackage{amssymb}
\usepackage{amsmath}
\usepackage{amsthm}
\usepackage{german}
\usepackage{psfrag}

\def\negthickspace{\!\!\!}

\newtheorem*{theorem}{Theorem}
\newtheorem*{corollary}{Corollary}
\newtheorem*{lemma}{Lemma}

\theoremstyle{definition}
\newtheorem*{remarks}{Remarks}

\newcommand{\nicefrac}[2]
{\leavevmode \kern.1em\raise.5ex\hbox{\the\scriptfont0 #1}
             \kern-.1em/\kern-.15em\lower.25ex
             \hbox{\the\scriptfont0 #2}}

\begin{document}

\begin{center}
{\Large{\bf A note on two-dimensional minimal surface graphs in $\mathbb R^n$}}\\[0.2cm]
{\Large{\bf and a theorem of Bernstein-Liouville type}}\\[0.8cm]
{\large{\sc Steffen Fr\"ohlich}}\\[1cm]
{\small\bf Abstract}\\[0.4cm]
\begin{minipage}[c][2cm][l]{12cm}
{\small Using Schauder's theory for linear elliptic partial differential equations in two independent variables and fundamental estimates for univalent mappings due to E. Heinz we establish an upper bound of the Gaussian curvature of two-dimensional minimal surface graphs in $\mathbb R^n.$ This leads us to a theorem of Bernstein-Liouville type.}
\end{minipage}
\end{center}
\vspace*{0.8cm}
{\small Acknowledgement: I am grateful for many discussions with Matthias Bergner.}\\[0.1cm]
\vspace{0.3cm}
\subsection{Non-parametric settings}
For natural $n\ge 3,$ let us consider a two-dimensional minimal surface graph
  $$X(x,y)=(x,y,\varphi_1(x,y),\ldots,\varphi_{n-2}(x,y)),\quad(x,y)\in\overline B_R\,,$$
with functions $\varphi_\Sigma\in C^2(\overline B_R,\mathbb R),$ $\Sigma=1,\ldots,n-2,$ on the topological closure $\overline B_R\subset\mathbb R^2$ of the open disc
  $$B_R:=\big\{(x,y)\in\mathbb R^2\,:\,x^2+y^2<R^2\big\},\quad R\in(0,+\infty).$$
With the minimal surface we associate the moving frame $\{X_x,X_y,N_1,\ldots,N_{n-2}\}$ consisting of the tangent vectors
  $$X_x=(1,0,\varphi_{1,x},\ldots,\varphi_{n-2,x}),\quad
    X_y=(0,1,\varphi_{1,y},\ldots,\varphi_{n-2,y}),$$
where the indices $x$ and $y$ denote the partial derivatives w.r.t. the coordinates, and the linearly independent unit normal vectors
\begin{equation}\label{4}
\begin{array}{rcl}
  N_1\negthickspace
  & := & \negthickspace\displaystyle
         \frac{1}{\sqrt{1+|\nabla\varphi_1|^2}}\,(-\varphi_{1,x},-\varphi_{1,y},1,0,\ldots,0,0), \\[0.6cm]
  \vdots & & \hspace*{3cm}\vdots \\[0.6cm]
  N_{n-2}\negthickspace
  & := & \negthickspace\displaystyle
         \frac{1}{\sqrt{1+|\nabla\varphi_{n-2}|^2}}\,(-\varphi_{n-2,x},-\varphi_{n-2,y},0,0,\ldots,0,1).
\end{array}
\end{equation}
\subsection{Introduction of conformal parameters}
Instead of the coordinates $x$ and $y,$ we prefer to consider the graph conformally parametrized (see \cite{Sauvigny_01}). Let $(u,v)\in\overline B,$ where $B:=B_1,$ then we get the immersion
  $$X(u,v)=(x^1(u,v),x^2(u,v),x^3(u,v),\ldots,x^n(u,v)),\quad(u,v)\in\overline B,$$
of regularity class $X\in C^2(B,\mathbb R^n)\cap C^0(\overline B,\mathbb R^n),$ while the following properties hold true:
\begin{itemize}
\item[(a)]
the plane mapping
  $$F^*(u,v):=(x^1(u,v),x^2(u,v))\,:\,\overline B\longrightarrow\overline B_R$$
is one-to-one, it has a positive Jacobian, it satisfies $F^*(0,0)=(0,0),$ and $F^*\,:\,\partial B\to\partial B_R$ is a topological mapping, where
  $$\partial B_R:=\big\{(x,y)\in\mathbb R^2\,:\,x^2+y^2=R^2\big\}$$
denotes the boundary of $B_R;$
\item[(b)]
the conformality relations
  $$|X_u|^2=W=|X_v|^2\,,\quad
    X_u\cdot X_v^t=0
    \quad\mbox{in}\ B$$
with the surface area element
  $$W:=\sqrt{h_{11}h_{22}-h_{12}^2}\,,\quad
    h_{ij}:=X_{u^i}\cdot X_{u^j}^t\quad\mbox{for}\ i,j=1,2$$
and $u^1\equiv u,$ $u^2\equiv v.$
\end{itemize}
Now, consider any differentiable unit normal vector field $N=N(u,v).$ W.r.t. its direction we define the coefficients of the second fundamental form by
  $$L_{N,ij}:=X_{u^iu^j}\cdot N^t\,,\quad i,j=1,2.$$
The mean curvature and the Gaussian curvature in direction of $N$ read
  $$H_N:=\frac{L_{N,11}h_{22}-2L_{N,12}h_{12}+L_{N,22}h_{11}}{2W^2}\,,\quad
    K_N:=\frac{L_{N,11}L_{N,22}-L_{N,12}^2}{W^2}\,.$$
In terms of the corresponding principle curvatures $\kappa_{N,1}$ and $\kappa_{N,2},$ we can write
  $$H_N=\frac{\kappa_{N,1}+\kappa_{N,2}}{2}\,,\quad
    K_N=\kappa_{N,1}\kappa_{N,2}\,.$$
\begin{lemma}
Let $X\in C^2(B,\mathbb R^n)\cap C^0(\overline B,\mathbb R^n)$ be a minimal immersion. Then there hold
  $$H_N(u,v)\equiv 0\quad\mbox{in}\ B$$
for all differentiable unit normal vector fields $N=N(u,v).$
\end{lemma}
\noindent
For the proof of this Lemma and further results on two-dimensional minimal immersions we refer to the fundamental paper \cite{Osserman}.
\subsection{The main result}
The main object of this note is the following
\begin{theorem}
Let the minimal surface graph $X(x,y)=(x,y,\varphi_1(x,y),\ldots,\varphi_{n-2}(x,y))$ with $\varphi_\Sigma\in C^2(\overline B_R,\mathbb R),$ $\Sigma=1,\ldots,n-2,$ be given. Then, for the principle curvatures $\kappa_{N,1}$ and $\kappa_{N,2}$ in direction of any unit normal vector field $N,$ there is a universal constant $\Theta\in(0,+\infty)$ such that in terms of conformal parameters $(u,v)\in\overline B$ it holds
  $$\kappa_{N,1}(0,0)^2+\kappa_{N,2}(0,0)^2
    \le\frac{1}{R^4}\,\Theta\,\|X\|_{C^0(B)}^2$$
with the Schauder norm
  $$\|X\|_{C^0(B)}:=\sup_{(x,y)\in B}|X(x,y)|.$$
\end{theorem}
\noindent
We conclude the following result of Bernstein-Liouville type:
\begin{corollary}
For large $R,$ let $\|X\|_{C^0(B_R)}\le\Omega R^\omega$ with a real constant $\Omega\in(0,+\infty)$ and $\omega\in[0,2).$ Then, if the minimal graph is complete, that is it is defined over the whole plane $\mathbb R^2,$ it represents an affine plane.
\end{corollary}
\noindent
First, we continue with some remarks.
\begin{remarks}
\begin{itemize}
\item[1.]
The corollary is sharp in the following sense: the holomorphic function $z\mapsto z^2,$ $z=x+iy,$ induces the minimal graph $(z,z^2)$ in $\mathbb R^4$ which is not linear.
\item[2.]
Graphs $(z,\Phi(z)),$ where $\Phi=\Phi(z)$ is any holomorphic function, are minimal in $\mathbb R^4.$ Thus, the above corollary can be read as a generalization of the well-known Liouville theorem.
\item[3.]
In \cite{Osserman}, using methods from complex function theory it was proved a curvature estimate and a theorem of Bernstein type for two-dimensional minimal immersions in $\mathbb R^n$ under the assumption that {\it all normals of such a surface omit any given neighborhood of some direction in space.}\\[0.1cm]
In \cite{Bergner_Froehlich}, Osserman's condition enabled us to establish upper curvature bounds for immersions of prescribed mean curvature fields. The assumptions of our theorem above are more restrictive but the proof of the theorem seems to be simpler.
\end{itemize}
\end{remarks}
\renewcommand{\proofname}{Proof of the Theorem}
\begin{proof}
\begin{itemize}
\item[1.]
Note that
  $$\begin{array}{lll}
      |K_N(0,0)|\negthickspace
      &  =  & \negthickspace\displaystyle
              \frac{|L_{N,11}(0,0)L_{N,22}(0,0)-L_{N,12}(0,0)^2|}{W(0,0)^2} \\[0.6cm]
      & \le & \negthickspace\displaystyle
              \frac{|X_{uu}(0,0)||X_{vv}(0,0)|+|X_{uv}(0,0)|^2}{W(0,0)^2}\,.
    \end{array}$$
Thus, we have to establish an upper bound for the second derivatives $X_{uu}(0,0)$ etc., and a lower bound for $W(0,0).$\\[1ex]
First, a minimal surface in conformal parameters is harmonic, that is (see \cite{Osserman})
  $$\triangle X(u,v)=0\quad\mbox{in}\ B.$$
Applying \cite{Gilbarg_Trudinger}, Theorem 4.6, there is a constant $C_1\in(0,+\infty)$ such that
  $$|X_{u^iu^j}(0,0)|\le C_1\|X\|_{C^0(B)},\quad i,j=1,2.$$
\item[2.]
Instead of the plane mapping $F^*=F^*(u,v)$ we consider the normalization
  $$F(u,v):=\frac{1}{R}\,F^*(u,v),\quad(u,v)\in\overline B.$$
Together with the properties of $F^*$ mentioned above in point (a), due to \cite{Sauvigny_02}, vol. II, chapter XII, Satz 1 there is a universal constant $C_2\in(0,+\infty)$ such that
  $$|\nabla F(0,0)|^2\ge C_2$$
(E. Heinz, 1952). Now, using the conformality relations we estimate as follows:
  $$\begin{array}{lll}
      2W(0,0)\negthickspace
      &  =  & \negthickspace\displaystyle
              |\nabla x^1(0,0)|^2+|\nabla x^2(0,0)|^2+|\nabla x^3(0,0)|^2+|\nabla x^4(0,0)|^2 \\[0.2cm]
      & \ge & \negthickspace\displaystyle
              |\nabla F^*(0,0)|^2
              \,=\,R^2|\nabla F(0,0)|^2
              \,\ge\,R^2 C_2\,.
    \end{array}$$
\item[3.]
We arrive at the estimate
  $$|K_N(0,0)|\le\frac{8C_1^2\|X\|_{C^0(B)}^2}{R^4C_2^2}\,.$$
Finally,
  $$\kappa_{N,1}^2+\kappa_{N,2}^2=2(2H_N^2-K_N)=2(-K_N)\ge 0\quad\mbox{for any}\ N$$
yields
  $$\kappa_{N,1}(0,0)^2+\kappa_{N,2}(0,0)^2
    \le\frac{1}{R^4}\,\frac{16C_1^2}{C_2^2}\,\|X\|_{C^0(B)}^2\,.$$
Setting $\Theta:=16C_1^2C_2^{-2}$ proves the theorem.
\end{itemize}
\end{proof}
\renewcommand{\proofname}{Proof of the Corollary}
\begin{proof}
We introduce conformal parameters $(u,v)\in\overline B.$
\begin{itemize}
\item[(i)]
Orthonormalization of the unit normal basis given in (\ref{4}) yields a unit normal frame $\{N_1(u,v),\ldots,N_{n-2}(u,v)\}$ such that
  $$N_\Sigma(u,v)\cdot N_\Omega(u,v)^t
    =\delta_{\Sigma\Omega}
     :=\left\{
         \begin{array}{l}
           1\quad\mbox{for}\ \Sigma=\Omega \\
           0\quad\mbox{for}\ \Sigma\not=\Omega
         \end{array}
       \right..$$
The Gaussian curvature $K$ (that is the non-trivial component of the Riemann curvature tensor) can be calculated by
  $$K=\sum_{\Sigma=1}^{n-2}K_\Sigma$$
with the Gauss curvatures $K_\Sigma$ in direction $N_\Sigma$ (see \cite{Brauner}, equation (6.3,17)). For our minimal surface graph there hold $H_\Sigma=0$ and $K_\Sigma\le 0$ for $\Sigma=1,\ldots,n-2.$
\item[(ii)]
For large $R,$ we have by assumption $\|X\|_{C^0(B)}\le\Omega R^\omega.$ From the curvature estimate above we conclude
  $$\kappa_{N,1}(0,0)^2+\kappa_{N,2}(0,0)^2
    \le\frac{1}{R^4}\,R^{2\omega}\Theta\Omega^2
    =\frac{1}{R^{4-2\omega}}\,\Theta\Omega^2$$
for any $N=N(u,v)$ of the orthonormal section $\{N_1,\ldots,N_{n-2}\}.$
Letting $R\to\infty$ implies $K_N(0,0)=0$ and therefore $K(0,0)=0$ for the Gaussian curvature. By \cite{Lawson}, chapter III, \S 5, Lemma 3, the graph is linear.
\end{itemize}
\end{proof}

\vspace*{1cm}
Steffen Fr\"ohlich\\
AG Differentialgeometrie und Geometrische Datenverarbeitung\\
TU Darmstadt\\
Schlo\ss{}gartenstra\ss{}e 7\\
D-64289 Darmstadt\\[0.2cm]
e-mail: sfroehlich@mathematik.tu-darmstadt.de

\end{document}